\documentclass[11pt]{article}
\usepackage{amssymb,latexsym,amsmath,amsbsy,amsthm,amsxtra,amsgen}
\numberwithin{equation}{section}

\newfont{\msbm}{msbm10 at 11pt}

\newcommand {\Z} {\mbox{\msbm Z}}

\def\be{\begin{equation}}
\def\ee{\end{equation}}
\def\ba{\begin{align}}
\def\ea{\end{align}}

\newtheorem{Theo}{Theorem}[section]
\newtheorem{Lemma}[Theo]{Lemma}

\newtheorem{Prop}[Theo]{Proposition}

\begin{document}

\title{Improving on bold play when the gambler is restricted}

\author{by Jason Schweinsberg\thanks{Supported by an NSF Postdoctoral Fellowship}}
\maketitle

\begin{abstract}
Suppose a gambler starts with a fortune in $(0,1)$ and wishes to
attain a fortune of $1$ by making a sequence of bets.  Assume that
whenever the gambler stakes the amount $s$, the gambler's fortune
increases by $s$ with probability $w$ and decreases by $s$ with
probability $1-w$, where $w < 1/2$.  Dubins and Savage showed that
the optimal strategy, which they called ``bold play'',
is always to bet $\min\{f, 1-f\}$, where
$f$ is the gambler's current fortune.  Here we consider the
problem in which the gambler may stake no more than $\ell$ at one
time.  We show that the bold strategy of always betting
$\min\{\ell, f, 1-f\}$ is not optimal if $\ell$ is
irrational, extending a result of Heath, Pruitt, and Sudderth.
\end{abstract}

\section{Introduction and background}
Suppose a gambler starts with a fortune in $(0,1)$ and wishes to
attain a fortune of $1$ by making a sequence of bets.  If the
gambler's current fortune is $f$, then the gambler may stake any
amount less than or equal to $f$.  The gambler wins the amount of the stake
with probability $w$ and loses the stake with probability $1-w$.
Following \cite{dusa76}, we refer to this game as red-and-black.
Clearly the gambler should never stake more than $1-f$, which is
enough to ensure that the gambler will reach the goal if the bet is won.
The strategy in which the gambler always stakes $\min\{f, 1-f\}$
is called bold play.

In \cite{dusa76}, Dubins and Savage developed a general theory for
gambling problems.  For red-and-black, they showed that if $0 < w < 1/2$,
which means that the game is subfair, then bold play is the optimal
strategy, in the sense that it maximizes the probability that the gambler
will eventually reach the goal.  Their proof is also given in
chapter 7 of \cite{bill95} and chapter 24 of \cite{fg97}.
See \cite{bak01} for some computations comparing the probability
that a gambler will reach the goal using bold play to the probability
that a gambler will reach the goal using other strategies.

This result has been extended in several ways.
Dubins and Savage \cite{dusa76} also considered
primitive casinos, in which the gambler loses the stake $s$
with probability $1-w$ and wins $[(1-r)/r]s$ with probability $w$,
where $0 < r < 1$.  Note that $r = 1/2$ is red-and-black.
They showed that bold play is optimal when the game is subfair, which in this
case means $w < r$.  Chen \cite{chen78} considered red-and-black
with inflation, in which the goal is not to reach $1$ but to
reach $(1 + \alpha)^n$ after $n$ bets
for some $n$.  He showed that bold play is optimal when $w \leq 1/2$.
A different extension is to incorporate a discount factor, so that
the gambler receives a utility of $\beta^n$, where
$0 < \beta \leq 1$, from reaching $1$ on the $n$th bet.  Klugman
\cite{klug77} showed that bold play is optimal for subfair
red-and-black with a discount factor.  However,
for some subfair primitive casinos, there exist discount factors
for which bold play is not optimal (see \cite{chen77} and
\cite{chen79}).  See also \cite{sec97} for a discussion of
the optimality of bold play in some two-person games.

Several authors have considered discrete versions of this problem in which
the gambler's initial fortune and the amount of each bet must be
integers and the gambler's goal is to attain a fortune of $n$.
An extensive discussion of discrete gambling problems
such as this can be found in \cite{masu96}.  Bold play remains optimal
when $w < 1/2$.  Ross \cite{ross74} showed that the timid strategy
of staking exactly $1$ each time is optimal in the superfair case
when $w > 1/2$.  See also \cite{ruth99} for an analysis of the
superfair case when the minimum bet is $2$.
Dubins \cite{dub98} showed, however, that if the win
probability is less than $1/2$ but is allowed to depend on the
gambler's fortune, then bold play need not be optimal.

Another direction of work concerns gambling problems in which there
is a limit to how much the gambler may bet.  The simplest
problem of this type involves red-and-black in which
the gambler may bet no more than $\ell \in (0, 1/2)$ at one time.  In this case,
we define bold play as the strategy in which the gambler whose
current fortune is $f$ always stakes $\min\{\ell, f, 1-f\}$.  Wilkins
\cite{wilk72} showed that if $w < 1/2$ and $\ell = 1/n$ for some
positive integer $n \geq 3$, then bold play maximizes the chance
that the gambler will reach the goal.  Chen \cite{chen76}
showed that bold play remains optimal when there is a discount factor
in addition to a limit on the stake of $1/n$.  In \cite{pesu85},
the optimality of bold play in continuous-time gambling problems
was established under rather general restrictions on the gambler.

However, Heath, Pruitt, and Sudderth \cite{hps72} obtained an
important negative result for discrete-time red-and-black.  They
showed that if the gambler can stake at most $\ell$, and if
$1/(n+1) < \ell < 1/n$ for some $n \geq 3$ or if $\ell$ is irrational
and $1/3 < \ell < 1/2$, then there exists $\epsilon > 0$ such that
if $0 < w < \epsilon$, then bold play is not optimal.  To see
heuristically why this is true, suppose $1/4 < \ell < 1/3$
and the gambler's initial fortune is $f = 1/2 - \delta$, where
$\delta$ is small.  If the gambler plays boldly and loses the
first bet, then the gambler's fortune after one bet will be
$1/2 - \ell - \delta$.  The gambler's fortune can at most double to
$1 - 2\ell - 2\delta$ after the second bet and therefore can be
at most $1 - 2\delta$ after two more wins.  However, if the
gambler first stakes $\ell - \delta$ and plays boldly thereafter,
then even with an initial loss, the gambler can reach the
goal by winning the next three bets.  Consequently, for sufficiently
small $\delta$, first betting $\ell - \delta$ makes the gambler more
likely to achieve the goal after winning three or fewer bets.
As $w \downarrow 0$, the probability that the gambler can win
four bets before going bankrupt gets very small relative to the
probability that the gambler wins three bets.  Therefore,
first betting $\ell - \delta$ is a better strategy than bold play
for sufficiently small $w$.

The purpose of the present paper is to extend this result
by showing that when $\ell$ is irrational, bold play fails to be optimal
for all $w < 1/2$, not just for very small $w$.  The case of
rational $\ell$ remains open except when $\ell = 1/n$ for some $n \geq 3$.
Note that when $\ell$ is rational and $1/3 < \ell < 1/2$, it is not even
known whether one can improve on bold play for very small $w$.

To state our result more precisely, define the function
$s: [0,1] \rightarrow [0,1]$ by
$s(f) = \min\{\ell, f, 1-f\}$.  We think of $s(f)$ as the bold stake for
a gambler whose fortune is $f$.  Denote by $X_k$ the gambler's
fortune after $k$ bets, when the gambler plays boldly.  Note that
$(X_k)_{k=0}^{\infty}$ is a Markov chain whose transition probabilities
are given by
\begin{align}
P(X_{k+1} = f + s(f)|X_k = f) &= w, \label{trans1} \\
P(X_{k+1} = f - s(f)|X_k = f) &= 1 - w. \label{trans2}
\end{align}
Define $Q(f) = P(X_k = 1 \mbox{ for some }k|X_0 = f)$, which is the
probability that a gambler who starts with a fortune of $f$ will
eventually reach the goal.  The following is our main result.

\begin{Theo}
Suppose $w < 1/2$ and $\ell$ is irrational.  Then, there exist
$f \in (0,1)$ and $\epsilon \in (0, s(f))$ such that
\begin{equation}
wQ(f + s(f) - \epsilon) + (1-w)Q(f - s(f) + \epsilon) > Q(f).
\label{maineq}
\end{equation}
\label{mainth}
\end{Theo}

If a gambler begins with a fortune of $f$ and stakes $s(f) - \epsilon$,
then the gambler's fortune after one bet will be $f + s(f) - \epsilon$
with probability $w$ and $f - s(f) + \epsilon$ with probability
$1 - w$.  Consequently, the left-hand side of (\ref{maineq}) is the
probability that the gambler will eventually reach the goal using the
strategy of first staking $s(f) - \epsilon$ and playing boldly
thereafter, while the right-hand side of (\ref{maineq}) is the
probability that the gambler will reach the goal using bold play.
Therefore, (\ref{maineq}) implies that the strategy of first
staking $s(f) - \epsilon$ and then playing boldly is superior to
bold play, and hence bold play is not optimal.

\section{Proof of Theorem \ref{mainth}}

\label{mainsec}

In this section, we prove Theorem \ref{mainth}.  The key to the
proof will be the following proposition.  Here, and throughout
the rest of the paper, all logarithms are assumed to be base 2.
That is, we write $\log n$ instead of $\log_2 n$.

\begin{Prop}
Let $S = \{f: P(X_k = 1 - \ell \mbox{ for some }k|X_0 = f) > 0\}$.
That is, $S$ is the set of all $f$ such that a gambler who
starts with a fortune of $f$ and plays boldly could have a
fortune of exactly $1-\ell$ after a finite number of bets.
\begin{enumerate}
\item Suppose $f \in S$.  Then there exists a constant
$C > 0$ such that if $0 < \epsilon < \ell$,
then $Q(f) - Q(f - \epsilon) \geq C (1 - w)^{-\log \epsilon}$.

\item Suppose $f \notin S$.  For all $C > 0$, there exists
$\delta > 0$ such that if $0 < \epsilon < \delta$, then
$Q(f) - Q(f - \epsilon) \leq C(1-w)^{-\log \epsilon}$.

\item If $\ell$ is irrational, then there exists $f \in (\ell, 1-\ell)$
such that $f - \ell \in S$ and $f + \ell \notin S$.
\end{enumerate}
\label{mainprop}
\end{Prop}

Proposition \ref{mainprop} implies that $Q(f) - Q(f - \epsilon)$
is larger when $f \in S$ than when $f \notin S$.
In other words, the difference between a fortune of $f$ and a
fortune of $f - \epsilon$ matters more to the gambler when
$f \in S$ than when $f \notin S$.
Part 3 of the proposition states that when $\ell$ is irrational,
we can find $f$ such that $f - s(f) \in S$ and $f + s(f) \notin S$.
We will show that if a gambler starts with a fortune
slightly below $f$, then it is better to make
slightly less than the bold stake so that the fortune
will not fall below $f - s(f)$ if the bet is lost.  This will
imply Theorem \ref{mainth}.

An important tool for the proof of Proposition \ref{mainprop}
is a coupling construction in which we follow two gamblers
simultaneously.  We present this construction in subsection 2.1.
We prove parts 1, 2, and 3 of Proposition \ref{mainprop} in
subsections 2.2, 2.3, and 2.4 respectively.  Then in subsection 2.5,
we show how Theorem \ref{mainth} follows from Proposition \ref{mainprop}.

\subsection{A coupling construction}
\label{css}

Throughout this and the next two subsections, we consider two Markov chains
$(X_k)_{k=0}^{\infty}$ and $(Y_k)_{k=0}^{\infty}$.  We define
$X_0 = f_1$ and $Y_0 = f_2$, where $f_1 \geq f_2$.
Both chains evolve with the transition
probabilities given by (\ref{trans1}) and (\ref{trans2}).
Consequently, we can think
of $X_k$ as the fortune after $k$ bets of a gambler whose initial
fortune is $f_1$, while $Y_k$ is the fortune after $k$ bets of a
gambler whose initial fortune is $f_2$.

We assume these sequences are coupled, so that both gamblers win
and lose the same bets.  To construct this coupling, we work with
the probability space $(\Omega, {\cal F}, P)$ defined as follows.
Let $\Omega = \{0,1\}^{\infty}$, and denote sequences in $\Omega$
by $\omega = (\omega_1, \omega_2, \ldots)$, so $\omega \rightarrow
\omega_i$ is the $i$th coordinate function.  Let ${\cal F}_0$ be
the trivial $\sigma$-field, and for positive integers $k$, let
${\cal F}_k$ be the
$\sigma$-field generated by the first $k$ coordinate functions.  Let
${\cal F} = \sigma({\cal F}_1, {\cal F}_2, \ldots)$ be the product
$\sigma$-field.  Let $P$ be the product probability measure with
the property that $P(\omega_i = 1) = w$ and $P(\omega_i = 0) = 1-w$
for all $i$.  We then say that the two gamblers win the $i$th bet
if $\omega_i = 1$ and lose the $i$th bet if $\omega_i = 0$.
In particular, for $k \geq 0$, we define $X_{k+1}(\omega) =
X_k(\omega) + s(X_k(\omega))$
and $Y_{k+1}(\omega) = Y_k(\omega) + s(Y_k(\omega))$
if $\omega_{k+1} = 1$, and we define
$X_{k+1}(\omega) = X_k(\omega) - s(X_k(\omega))$
and $Y_{k+1}(\omega) = Y_k(\omega) - s(Y_k(\omega))$ if
$\omega_{k+1} = 0$.

We now make some remarks pertaining to this construction.

\begin{enumerate}
\item Since $s(f) = \min\{\ell, f, 1-f\}$, we see that if
$f \leq g$, then $|s(g) - s(f)| \leq g - f$.  Therefore,
$f + s(f) \leq g + s(g)$ and $f - s(f) \leq g - s(g)$.  It then
follows by induction and the construction of the sequences
$(X_k)_{k=0}^{\infty}$ and $(Y_k)_{k=0}^{\infty}$ that
$X_k \geq Y_k$ for all $k$.  Likewise, the fact that
$|s(g) - s(f)| \leq g - f$ implies that
$X_k - Y_k \leq 2^k(f_1 - f_2)$ for all $k$.

\item The fact that $X_k \geq Y_k$ for all $k$ means that if
$Y_k = 1$ then $X_k = 1$.
Since $Q(f_1) = P(X_k = 1 \mbox{ for some }k)$ and
$Q(f_2) = P(Y_k = 1 \mbox{ for some }k)$, it follows that
$Q(f_1) \geq Q(f_2)$.  That is, the function
$f \mapsto Q(f)$ is nondecreasing.

\item Note that
\begin{align}
E[X_{k+1}|{\cal F}_k] &= w(X_k + s(X_k)) + (1-w)(X_k - s(X_k)) \nonumber \\
&= X_k + (2w - 1)s(X_k) \leq X_k, \nonumber
\end{align}
where the last inequality holds because $w < 1/2$.  Therefore,
$(X_k)_{k=0}^{\infty}$ is a supermartingale with respect to
$({\cal F}_k)_{k=0}^{\infty}$.  By the same argument,
$(Y_k)_{k=0}^{\infty}$ is a supermartingale with respect to
$({\cal F}_k)_{k=0}^{\infty}$.  By the Martingale Convergence Theorem
(see chapter 4 of \cite{durr96}), there exist random variables
$L_1$ and $L_2$ such that $X_k \rightarrow L_1$ a.s. and
$Y_k \rightarrow L_2$ a.s. as $k \rightarrow \infty$.  
If $0 < \epsilon < \ell$, then $s(f) > \epsilon$ for
$f \in [\epsilon, 1 - \epsilon]$.  It follows that
$L_1$ and $L_2$ must be $\{0,1\}$-valued random variables.
Furthermore, it is easy to see that $X_k = 1$ for sufficiently large
$k$ on $\{L_1 = 1\}$ and $Y_k = 1$ for sufficiently large $k$
on $\{L_2 = 1\}$.  Thus, $Q(f_1) = P(L_1 = 1)$ and
$Q(f_2) = P(L_2 = 1)$, from which it follows that
$Q(f_1) - Q(f_2) = P(L_1 = 1 \mbox{ and }L_2 = 0)$.
\end{enumerate}

\subsection{Proof of part 1 of Proposition \ref{mainprop}}

We begin with the following Lemma, in which we compute the gambler's
probability of reaching the goal starting from a
sequence of fortunes approaching $1$.

\begin{Lemma}
For all $n \geq 0$, we have $Q(1 - 2^{-n}\ell) =
1 - (1-w)^n (1 - Q(1-\ell))$.
\label{Q2nlem}
\end{Lemma}

\begin{proof}
The statement is obvious when $n = 0$.  Suppose the result
holds for some $n \geq 0$.  Since $s(1 - 2^{-(n+1)}\ell) = 2^{-(n+1)}\ell$,
a gambler whose fortune is $1 - 2^{-(n+1)}\ell$ will, after the next bet,
have a fortune of $1$ with probability $w$ and a fortune of 
$1 - 2^{-n} \ell$ with probability $1 - w$.  Thus, by the Markov property,
\begin{align}
Q(1 - 2^{-(n+1)}\ell) &= w + (1-w)Q(1 - 2^{-n}\ell) \nonumber \\
&= w + (1-w)(1 - (1-w)^n (1 - Q(1-\ell))) \nonumber \\
&= 1 - (1-w)^{n+1}(1 - Q(1-\ell)). \nonumber
\end{align}
The lemma now follows by induction on $n$.
\end{proof}

\bigskip
\noindent {\it Proof of part 1 of Proposition \ref{mainprop}.}  
Let $f_1 = f$ and $f_2 = f - \epsilon$, where $0 < \epsilon < \ell$.
Since $f \in S$, there exists a positive integer $k$ such that if
$B$ denotes the event that $X_k = 1 - \ell$ and $X_{k+1} = 1$, then
$P(B) > 0$.  Note that for $0 \leq j < k$, we have
$X_{j+1} - Y_{j+1} \geq X_j - Y_j$ unless either
$X_j > 1 - \ell$ and $X_{j+1} = 1$, or $Y_j < \ell$ and $Y_{j+1} = 0$.
Therefore, if $B$ occurs, then $Y_k \leq 1 - \ell - \epsilon$ and
thus $Y_{k+1} \leq 1 - \epsilon$.  Combining this observation with
remarks 2 and 3 in subsection \ref{css}, we get
$$Q(f) - Q(f - \epsilon) = P(L_1 = 1 \mbox{ and } L_2 = 0) \geq
P(B)P(L_2 = 0|B) \geq P(B)(1 - Q(1 - \epsilon)).$$
Choose a nonnegative integer
$n$ such that $2^{-(n+1)}\ell < \epsilon \leq 2^{-n}\ell$,
which implies that $n \leq \log \ell - \log \epsilon$.  By Lemma \ref{Q2nlem},
\begin{align}
Q(1 - \epsilon) &\leq Q(1 - 2^{-(n+1)}\ell) = 1 - (1 - w)^{n+1}
(1 - Q(1-\ell)) \nonumber \\
&\leq 1 - (1 - w)^{1 + \log \ell - \log \epsilon}(1 - Q(1-\ell)). \nonumber
\end{align}
Thus, $Q(f) - Q(f - \epsilon) \geq C(1 - w)^{-\log \epsilon}$,
where $C = P(B)(1-w)^{1 + \log \ell}(1 - Q(1-\ell))$.

\subsection{Proof of part 2 of Proposition \ref{mainprop}}

Our next step is to prove part 2 of Proposition \ref{mainprop},
which gives an upper bound for $Q(f) - Q(f - \epsilon)$
when $f \notin S$.  We will compare the sequences
$(X_k)_{k=0}^{\infty}$ and $(Y_k)_{k=0}^{\infty}$ when
$f_1 = f$ and $f_2 = f - \epsilon$.  Although
$(X_k - Y_k)_{k=0}^{\infty}$ is not a supermartingale,
we will be able to construct a supermartingale
by considering the differences between the gamblers' fortunes at
a sequence of stopping times.  It will then follow that the gamblers'
fortunes stay close enough together for us to obtain the desired
upper bound on $Q(f) - Q(f - \epsilon)$ when $f \notin S$.

Given $f$ and $f^*$ such that
$0 \leq f^* \leq f \leq 1$, define
\begin{displaymath}
h(f, f^*) = \left\{
\begin{array}{ll} 1 & \mbox{ if }f = f^* \\
(s(f) - s(f^*))/(f - f^*) & \mbox{ otherwise. }
\end{array} \right.
\end{displaymath}
Note that $-1 \leq h(f, f^*) \leq 1$ for all $f$ and $f^*$.
If $\ell \leq f^* \leq f \leq 1-\ell$, then $s(f) = s(f^*) = \ell$, which
means $h(f, f^*) = 0$.  If $f^* \geq \ell$ and $f \geq 1-\ell$, then
$h(f, f^*) \leq 0$.  If $f^* \leq \ell$ and $f \leq 1- \ell$, then
$h(f, f^*) \geq 0$.  Also, recall that $\omega_k = 1$ if the
gamblers win the $k$th bet, and $\omega_k = 0$ if the gamblers lose
the $k$th bet.  We have
\begin{equation}
X_{k+1}(\omega) - Y_{k+1}(\omega) = \left\{
\begin{array}{ll} (1 + h(X_k(\omega), Y_k(\omega)))(X_k(\omega) - Y_k(\omega))
& \mbox{ if }\omega_{k+1} = 1 \\
(1 - h(X_k(\omega), Y_k(\omega)))(X_k(\omega) - Y_k(\omega)) &
\mbox{ if }\omega_{k+1} = 0.
\end{array} \right.
\label{hdis}
\end{equation}
Define
\begin{equation}
W_k = (1-w)^{-\log (X_k - Y_k)} = (X_k - Y_k)^{-\log(1-w)}.
\label{defW}
\end{equation}
By (\ref{hdis}), we have
\begin{align}
E[W_{k+1}|{\cal F}_k]
&= w(1 + h(X_k, Y_k))^{-\log(1-w)}W_k + (1-w)(1 - h(X_k, Y_k))^
{-\log(1-w)}W_k \nonumber \\
&= g(h(X_k, Y_k))W_k,
\label{gyeq}
\end{align}
where
\begin{displaymath}
g(x) = w(1 + x)^{-\log(1-w)} + (1-w)(1 - x)^{-\log(1-w)}
\end{displaymath}
for $-1 \leq x \leq 1$.
Note that
\begin{displaymath}
g'(x) = (-\log (1-w)) \big( w(1 + x)^{- \log(1-w) - 1} -
(1-w)(1 - x)^{- \log (1-w) - 1} \big).
\end{displaymath}
Suppose $0 < x < 1$.
Since $0 < -\log (1-w) < 1$, we have $(1 + x)^{-\log(1-w) - 1} < 1$
and $(1-x)^{-\log(1-w) - 1} > 1$.  Therefore,
\begin{displaymath}
g'(x) \leq (-\log(1-w))(w - (1-w)) < 0.
\end{displaymath}
Since $g(0) = 1$, it follows that $0 < g(x) < 1$ for $x \in (0,1]$.

We now introduce four lemmas that will help us to define a
supermartingale.

\begin{Lemma}
Suppose $1-\ell \leq f_2 \leq f_1 \leq 1$.  Then $E[W_1] = W_0$.
\label{LemA}
\end{Lemma}

\begin{proof}
We have $s(f_1) = 1 - f_1$ and $s(f_2) = 1 - f_2$.
Therefore, $h(f_1, f_2) = -1$.  Since $g(-1) = 1$, it follows
from (\ref{gyeq}) that $E[W_1] = W_0$.
\end{proof}

\begin{Lemma}
Suppose $f_2 \leq f_1 < 1-\ell$.  Define a stopping time $R$ as follows.
If $h(f_1, f_2) \geq 1/2$, define $R = 0$.  If $h(f_1, f_2) < 1/2$,
then let $R(\omega) = \inf\{j: \omega_j = 1 \mbox{ or }
h(X_j(\omega), Y_j(\omega)) \geq 1/2\}$.
Let $L = \lfloor 1 + (1 - 2\ell)/\ell \rfloor$.
Then $R \leq L$ and $E[W_R] \leq W_0$.
\label{Rlem}
\end{Lemma}

\begin{proof}
Proceeding by contradiction, suppose $R(\omega) > L$
for some $\omega$.  Then the gamblers
must lose the first $L$ bets.  However, by the definition of $L$, any gambler
who starts with a fortune of at most $1-\ell$ and then loses $L$ consecutive
bets has a fortune of at most $\ell$.  Therefore, there exists $j \leq L$
such that $0 < X_j \leq \ell$.  Since $X_j \geq Y_j$, it follows that
$Y_j \leq \ell$, and thus $s(X_j) = X_j$ and $s(Y_j) = Y_j$.  However,
this means that $h(X_j, Y_j) = 1$, and thus $R \leq j$, a contradiction.
Hence, $R \leq L$.

For $j < R$, we have $Y_j \leq X_j < 1-\ell$, and therefore
$0 \leq h(X_j, Y_j) \leq 1$.  Since $g(x) \leq 1$ for $x \in [0,1]$, we have,
with the aid of (\ref{gyeq}),
\begin{align}
E[W_{(j+1) \wedge R}|{\cal F}_j] &= W_{j \wedge R} 1_{\{R \leq j\}} +
E[W_{j+1}|{\cal F}_j]1_{\{R > j\}} \nonumber \\
&= W_{j \wedge R} 1_{\{R \leq j\}} + g(h(X_j, Y_j)) W_j 1_{\{R > j\}} \leq
W_{j \wedge R}. \nonumber
\end{align}
Therefore, $(W_{j \wedge R})_{j=0}^{\infty}$ is a supermartingale
with respect to $({\cal F}_j)_{j=0}^{\infty}$.
Note that $0 \leq W_{j \wedge R} \leq 1$ for all $j$, so the
Optional Stopping Theorem (see chapter 4 of \cite{durr96}) gives
$E[W_R] \leq W_0$.
\end{proof}

\begin{Lemma}
Suppose $f_2 \leq f_1 < 1-\ell$.  Let
$$T(\omega) = \inf\{j \geq 1: \omega_j = 1 \mbox{ or }h(X_{j-1}(\omega),
Y_{j-1}(\omega)) \geq 1/2\}.$$
Let $L = \lfloor 1 + (1 - 2\ell)/\ell \rfloor$ as in Lemma \ref{Rlem}.  Then
$T \leq L + 1$ and $E[W_T] \leq \alpha W_0$, where
$$\alpha = 1 - (1 - g(1/2))(1-w)^{2L}.$$
\label{LemB}
\end{Lemma}

\begin{proof}
Define the stopping time $R$ as in Lemma \ref{Rlem}.
Then $T = R$ if and only if the gamblers win the $R$th bet;
otherwise, $T = R+1$.  Clearly $T \leq L+1$ by Lemma \ref{Rlem}.
Let $A$ be the event that the gamblers win the $R$th bet.  Then
$$E[W_T] = E[E[W_T|{\cal F}_R]] = E[W_R 1_A + E[W_{R+1}|{\cal F}_R]1_{A^c}].$$
By the strong Markov property and (\ref{gyeq}), $E[W_{R+1}|{\cal F}_R] =
g(h(X_R, Y_R))W_R$.  If the gamblers lose the $R$th bet, then
$h(X_R, Y_R) \geq 1/2$.  Therefore, since $g$ is decreasing on $[0,1]$,
$$E[W_{R+1}|{\cal F}_R]1_{A^c} \leq g(1/2)W_R 1_{A^c}.$$  Thus,
\begin{equation}
E[W_T] \leq E[W_R 1_A + g(1/2) W_R 1_{A^c}] =
E[W_R - (1 - g(1/2))W_R 1_{A^c}].
\label{ew1}
\end{equation}

If $A^c$ occurs, then the gamblers lose the first $R$ bets, and
$h(X_j, Y_j) < 1/2$ for all $j < R$.  If $h(X_j, Y_j) < 1/2$ and the
gamblers lose the $(j+1)$st bet, then $X_{j+1} - Y_{j+1} \geq
(X_j - Y_j)/2$.  Thus, on $A^c$ we have $X_R - Y_R \geq
2^{-R}(f_1 - f_2) \geq 2^{-L}(f_1 - f_2)$.  Thus,
\begin{align}
E[W_R 1_{A^c}] &\geq E[(2^{-L}(f_1 - f_2))^{-\log (1-w)} 1_{A^c}] \nonumber \\
&= 2^{L \log (1-w)} W_0 P(A^c) = P(A^c)(1-w)^L W_0.
\label{ew2}
\end{align}
Since $A^c$ occurs whenever the gamblers lose
the first $L$ bets, we have $P(A^c) \geq (1-w)^L$.  Thus,
since $E[W_R] \leq W_0$ by Lemma \ref{Rlem}, 
combining (\ref{ew1}) and (\ref{ew2}) gives
$$E[W_T] \leq E[W_R] - (1 - g(1/2))(1 - w)^{2L}W_0 \leq
(1 - (1 - g(1/2))(1-w)^{2L})W_0 = \alpha W_0,$$
which completes the proof.
\end{proof}

\begin{Lemma}
Suppose $f_2 < 1-\ell$ and $1-\ell \leq f_1 < 1 - \ell/2$.  
Define the stopping time $T$ by
\begin{displaymath}
T(\omega) = \left\{
\begin{array}{ll} 1 & \mbox{ if }\omega_1 = 1 \\
\inf\{j \geq 2: \omega_j = 1 \mbox{ or }h(X_{j-1}(\omega), Y_{j-1}(\omega))
\geq 1/2 \}
& \mbox{ otherwise. }
\end{array} \right.
\end{displaymath}
Let $N(\omega) = 1 - \alpha$ if $\omega_1 = 1$, and let $N(\omega) = 1$ if
$\omega_1 = 0$.  Then $T \leq L + 2$ and $E[NW_T] \leq W_0$.
\label{LemC}
\end{Lemma}

\begin{proof}
Let $A$ be the event that that the gamblers win the first bet,
which means $\omega_1 = 1$.  We have
$$E[NW_T] = E[E[NW_T|{\cal F}_1]] =
E[(1 - \alpha)W_1 1_A + E[W_T|{\cal F}_1] 1_{A^c}].$$
If the gamblers lose the first bet, then $X_1 = f_1 - s(f_1) =
2f_1 - 1 < 1-\ell$.  Therefore, Lemma \ref{LemB} and the Markov property give
$T \leq L + 2$ and $E[W_T|{\cal F}_1]1_{A^c} \leq \alpha W_1 1_{A^c}$.  Thus,
\begin{align}
E[NW_T] &\leq E[(1 - \alpha)W_1 1_A +
\alpha W_1 1_{A^c}] \nonumber \\
&= w(1 - \alpha)[(1 + h(f_1, f_2))(f_1 - f_2)]^{-\log(1-w)} \nonumber \\
&\hspace{.3in} +
(1-w)\alpha[(1 - h(f_1, f_2))(f_1 - f_2)]^{-\log (1-w)} \nonumber \\
&= \big[w(1 - \alpha) (1 + h(f_1, f_2))^{-\log(1-w)} +
(1-w)\alpha (1 - h(f_1, f_2))^{-\log(1-w)} \big] W_0 \nonumber \\
&\leq \big[w(1-\alpha)2^{-\log(1-w)} + (1-w) \alpha 2^{-\log(1-w)} \big]W_0
\nonumber \\
&= \bigg( \frac{w}{1-w} (1-\alpha) + \alpha \bigg) W_0 \leq W_0, \nonumber
\end{align}
which completes the proof.
\end{proof}

\bigskip
By combining Lemmas \ref{LemA}, \ref{LemB}, and \ref{LemC}, we can
obtain Proposition \ref{martprop}, in which we construct the
supermartingale needed to prove part 2 of Proposition \ref{mainprop}.
We first inductively define a sequence of stopping times
$(T_k)_{k=0}^{\infty}$.  Let $T_0 = 0$.  Given $T_k$, we define
$T_{k+1}$ according to the following rules:
\begin{enumerate}
\item If $Y_{T_k}(\omega) \geq 1-\ell$, then let
$T_{k+1}(\omega) = T_k(\omega) + 1$.

\item If $X_{T_k}(\omega) < 1-\ell$, let $$T_{k+1}(\omega) =
\inf\{j \geq T_k(\omega) + 1:
\omega_j = 1 \mbox{ or }h(X_{j-1}(\omega), Y_{j-1}(\omega)) \geq 1/2\}.$$

\item Suppose $Y_{T_k}(\omega) < 1-\ell$ and $1-\ell \leq X_{T_k}(\omega) < 1 - \ell/2$.
If $\omega_{T_k(\omega) + 1} = 1$, meaning the gamblers win the
$(T_k + 1)$st bet, then let $T_{k+1}(\omega) = T_k(\omega) + 1$.
Otherwise, let $T_{k+1}(\omega) = \inf\{j \geq T_k(\omega) + 2:
\omega_j = 1 \mbox{ or } h(X_{j-1}(\omega), Y_{j-1}(\omega)) \geq 1/2\}$.

\item If $Y_{T_k} < 1-\ell$ and $X_{T_k} \geq 1 - \ell/2$, then let
$T_{k+1} = T_k$.
\end{enumerate}

\begin{Prop}
Define the sequence of stopping times $(T_k)_{k=0}^{\infty}$ as above.
For $k \geq 0$, let $B_k = \#\{j \in \{0, 1, \ldots, k-1\}:
X_{T_j} < 1-\ell \}$, where $\#S$ denotes the cardinality of the
set $S$.  Let $N_k = 1 - \alpha$ on the event that,
for some $j \leq T_k$, we have $Y_{j-1} < 1- \ell \leq X_{j-1}$ and
the gamblers win the $j$th bet.  Otherwise, let $N_k = 1$.  Define
$Z_k = \alpha^{-B_k} N_k W_{T_k}$.  Then
$(Z_k)_{k=0}^{\infty}$ is a supermartingale with respect to the
filtration $({\cal F}_{T_k})_{k=0}^{\infty}$.  Furthermore,
$T_{k+1} \leq T_k + L + 2$ for all $k \geq 0$.
\label{martprop}
\end{Prop}

\begin{proof}
Let $A_{1,k}$ be the event that $Y_{T_k} \geq 1 - \ell$, and
let $A_{2,k}$ be the event that $X_{T_k} < 1- \ell$.  Let $A_{3,k}$ be the
event that $Y_{T_k} < 1-\ell$ and $1-\ell \leq X_{T_k} < 1 - \ell/2$. 
Let $A_{4,k}$ be the event that $Y_{T_k} < 1-\ell$ and $X_{T_k} \geq 1 - \ell/2$.
Note that for all $k$, exactly one of these four events occurs.
We consider the four cases separately.

First, suppose $A_{1,k}$ occurs.  Then $X_{T_k} \geq 1- \ell$, so
$B_{k+1} = B_k$.  Also, note that $T_{k+1} = T_k + 1$ and
$Y_{T_k} \geq 1-\ell$, so $N_{k+1} = N_k$.  Therefore, by Lemma 
\ref{LemA} and the strong Markov property,
\begin{align}
E[Z_{k+1} 1_{A_{1,k}}|{\cal F}_{T_k}] &=
\alpha^{-B_k} N_k E[W_{T_{k+1}}|{\cal F}_{T_k}] 1_{A_{1,k}} \nonumber \\
&= \alpha^{-B_k} N_k E[W_{T_k + 1}|{\cal F}_{T_k}] 1_{A_{1,k}} \nonumber \\
&= \alpha^{-B_k} N_k W_{T_k} 1_{A_{1,k}} = Z_k 1_{A_{1,k}}.
\label{qq1}
\end{align}

Next, suppose $A_{2,k}$ occurs.  Then $X_{T_k} < 1-\ell$, so
$B_{k+1} = B_k + 1$.  The gamblers lose bets $T_k + 1, \ldots,
T_{k+1} - 1$, so $X_j < 1-\ell$ for $T_k \leq j \leq T_{k+1} - 1$.
Therefore, $N_{k+1} = N_k$.  By Lemma \ref{LemB} and the strong
Markov property,
\begin{align}
E[Z_{k+1} 1_{A_{2,k}}|{\cal F}_{T_k}] &=
\alpha^{-(B_k + 1)} N_k E[W_{T_{k+1}}|{\cal F}_{T_k}] 1_{A_{2,k}}
\nonumber \\
&\leq \alpha^{-(B_k + 1)} N_k (\alpha W_{T_k}) 1_{A_{2,k}} =
\alpha^{-B_k} N_k W_{T_k} 1_{A_{2,k}} = Z_k 1_{A_{2,k}}.
\label{qq2}
\end{align}

Suppose $A_{3,k}$ occurs.  Then $X_{T_k} \geq 1- \ell$, so
$B_{k+1} = B_k$.  Since $X_{T_k} < 1$, we have $N_k = 1$.
If the gamblers win the $(T_k + 1)$st bet, then $N_{k+1} = 1 - \alpha$.
Otherwise, $X_{T_k + 1} = X_{T_k} - s(X_{T_k}) = 2X_{T_k} - 1 < 1- \ell$
and the gamblers lose bets $T_k + 2, \ldots, T_{k+1} - 1$, so
$N_{k+1} = 1$.  By Lemma \ref{LemC} and the strong Markov property,
\begin{equation}
E[Z_{k+1} 1_{A_{3,k}}|{\cal F}_{T_k}] = \alpha^{-B_k}
E[N_{k+1} W_{T_{k+1}}|{\cal F}_{T_k}] 1_{A_{3,k}} \leq
\alpha^{-B_k} W_{T_k} 1_{A_{3,k}} = Z_k 1_{A_{3,k}}.
\label{qq3}
\end{equation}

Finally, if $A_{4,k}$ occurs, then $T_{k+1} = T_k$.  Therefore,
$N_{k+1} = N_k$ and $W_{T_{k+1}} = W_{T_k}$.  Since $X_{T_k} > 1- \ell$, we also
have $B_{k+1} = B_k$ and thus $Z_{k+1} = Z_k$.  Therefore,
$E[Z_{k+1} 1_{A_{4,k}}|{\cal F}_{T_k}] = Z_k 1_{A_{4,k}}$.  
This fact, combined with equations (\ref{qq1}), (\ref{qq2}),
and (\ref{qq3}), gives $E[Z_{k+1}|{\cal F}_{T_k}] \leq Z_k$.
Hence, $(Z_k)_{k=0}^{\infty}$ is a supermartingale with
respect to the filtration $({\cal F}_{T_k})_{k=0}^{\infty}$.

To complete the proof, note that clearly $T_{k+1} \leq T_k + L + 2$
if $A_{1,k}$ or $A_{4,k}$ occurs.  The strong Markov property
combined with Lemmas \ref{LemB} and \ref{LemC} implies that
$T_{k+1} \leq T_k + L + 2$ if $A_{2,k}$ or $A_{3,k}$ occurs.
\end{proof}

\bigskip
We now use Proposition \ref{martprop} to establish an upper bound on
$Q(f) - Q(f - \epsilon)$ when $f \notin S$.  We will need one more lemma.

\begin{Lemma}
Fix $f \notin S$, and let $N$ be a positive integer.  Then there exists
a positive integer $M$ and a positive real number $\delta$ such
that if $f_1 = f$ and $f_2 = f - \epsilon$ where $0 < \epsilon < \delta$,
then the following hold:
\begin{enumerate}
\item If $X_k \geq 1-\ell$ and $Y_k < 1- \ell$, then $k > M(L + 2)$, where
$L = \lfloor 1 + (1 - 2\ell)/ \ell \rfloor$.

\item Let $D_k = \# \{j \in \{0, 1, \ldots, k-1\}:
X_{T_j} < 1- \ell \mbox{ or }X_{T_j} = 1\}$.  Then $D_M \geq N$.
\end{enumerate}
\label{endlem}
\end{Lemma}

\begin{proof}
Let $R_k'$ be the set of all possible values of $X_k$, and
let $R_k = \bigcup_{j=0}^k R_j'$.  Note that $R_k$ is a finite set
because there are only $2^k$ possible outcomes for the first $k$ bets.
For all $g \in [0, 1)$, let $v(g)$ be the number of consecutive
bets that a gambler whose fortune is $g$ must lose for the fortune to
drop below $1- \ell$.  That is, $v(g) = 0$ when $0 \leq g < 1- \ell$ and,
for positive integers $k$, $v(g) = k$ when
$1 - 2^{-k+1} \ell \leq g < 1 - 2^{-k} \ell$.  Let $V_k = \max\{v(g): g \in R_k\}$.
Let $M_0 = 0$.  For $i \geq 0$, let $M_{i+1} = M_i + V_{M_i(L+2)} + 1$.
Let $M = M_N$.  Choose $\theta > 0$ small enough that 
$R_{M(L+2)} \cap (1- \ell, 1- \ell + \theta) = \emptyset$.  Let
$\delta = 2^{-M(L+2)}\theta$.  We will show that the two conditions of
the lemma are satisfied for these choices of $M$ and $\delta$.

Suppose $k \leq M(L+2)$.  If $X_k \geq 1- \ell$, then $X_k \geq 1 - \ell + \theta$,
since $f \notin S$ and $R_k' \cap (1- \ell, 1- \ell + \theta) = \emptyset$.
By the first remark in subsection \ref{css}, we have
$X_k - Y_k \leq 2^k \epsilon < 2^{M(L+2)} \delta = \theta$.
Therefore, $Y_k \geq 1- \ell$.  This proves the first part of the lemma.

To prove the second part, we claim that for $i = 0, 1, \ldots, N-1$, we have
\begin{equation}
D_{M_i + V_{M_i(L+2)} + 1} \geq D_{M_i} + 1.
\label{Deq}
\end{equation}
To see how (\ref{Deq}) implies the second part of the lemma, first note
that $D_{M_0} = D_0 = 0$.  Suppose $D_{M_i} \geq i$ for some
$i \geq 0$.  Then $D_{M_{i+1}} = D_{M_i + V_{M_i(L+2)} + 1} \geq
D_{M_i} + 1 \geq i+1$ by (\ref{Deq}).
Hence, by induction, (\ref{Deq}) implies that
$D_M \geq N$.  Thus, we need only to prove (\ref{Deq}).  
First, suppose either $X_{T_{M_i}} < 1- \ell$ or $X_{T_{M_i}} = 1$.
Then $D_{M_i + 1} = D_{M_i} + 1$.  Since $(D_i)_{i=0}^{\infty}$
is a nonincreasing sequence and $V_{M_i(L+2)} \geq 0$, we have (\ref{Deq}).

Thus, it remains only to prove (\ref{Deq}) when 
$1- \ell \leq X_{T_{M_i}} < 1$.  Write $v$ for $v(X_{T_{M_i}})$.
Note that $v \leq V_{T_{M_i}}$, and $T_{M_i} \leq M_i(L+2)$ by
Proposition \ref{martprop}.  Therefore
\begin{align}
T_{M_i} + v &\leq T_{M_i} + V_{T_{M_i}} \leq M_i(L+2) + V_{M_i(L+2)}
\nonumber \\
&\leq (M_i + V_{M_i(L+2)})(L+2) \leq M(L+2).
\label{tM}
\end{align}
We now consider two cases.  First, suppose the gamblers lose the bets
$T_{M_i} + 1, \ldots, T_{M_i} + v$.  Then 
$X_j \geq 1- \ell$ for $T_{M_i} \leq j \leq T_{M_i + v - 1}$ and
$X_{T_{M_i} + v} < 1- \ell$.  Also, by (\ref{tM}) and the first part
of the lemma, we have $Y_j \geq 1- \ell$ for $T_{M_i} \leq j \leq T_{M_i + v - 1}$.
Therefore, by the definition of the sequence $(T_j)_{j=0}^{\infty}$,
we have $T_{M_i + k} = T_{M_i} + k$ for $1 \leq k \leq v$.
It follows that $X_{T_{M_i + v}} < 1- \ell$, which means
$D_{M_i + v + 1} = D_{M_i} + 1$.  Thus,
$D_{M_i + V_{M_i(L + 2)} + 1} \geq D_{M_i + v + 1} = D_{M_i} + 1$,
which is (\ref{Deq}).  Finally, we consider the case in which,
for some $j \in \{1, \ldots, v\}$, the gamblers lose the bets
$T_{M_i} + 1, \ldots, T_{M_i} + j - 1$ but win the bet $T_{M_i} + j$.
Then, $X_{T_{M_i} + j} = X_{T_{M_i + j}} = 1$ and
$D_{M_i + j + 1} = D_{M_i} + 1$.  Hence,
$D_{M_i + V_{M_i(L + 2)} + 1} \geq D_{M_i} + 1$, which is (\ref{Deq}).
\end{proof}

\bigskip
\noindent {\it Proof of part 2 of Proposition \ref{mainprop}.}
Fix $C > 0$, and fix $f \notin S$.
Since $0 < \alpha < 1$, there exists a positive integer
$N$ such that  $$\alpha^{-N} (1 - \alpha) (1-w)^{-\log (\ell/2)}
\geq C^{-1}.$$  Define $M$ and $\delta$ as in Lemma \ref{endlem}.
Fix $\epsilon \in (0, \delta)$.  Define $(T_k)_{k=0}^{\infty}$ and
$(Z_k)_{k=0}^{\infty}$ as in Proposition \ref{martprop}, with
$f_1 = f$ and $f_2 = f - \epsilon$.

By Remark 3 in subsection \ref{css}, there exist random variables 
$L_1$ and $L_2$ such that $X_k \rightarrow L_1$ a.s. and
$Y_k \rightarrow L_2$ a.s. as $k \rightarrow \infty$, and
$Q(f) - Q(f - \epsilon) = P(L_1 = 1 \mbox{ and }L_2 = 0)$.
Let $A$ be the event that $L_1 = 1$ and $L_2 = 0$.  Then there
is an integer-valued random variable $K$ such that, on the event $A$,
we have $X_{T_K} \geq 1 - \ell/2$ and $Y_{T_K} < 1- \ell$.
By part 1 of Lemma \ref{endlem}, on the event $A$, we have $T_K > M(L+2)$
and thus $K \geq M$.  It also follows from part 1 of Lemma \ref{endlem}
that if $X_{T_j} = 1$ for $j \leq M$, then $Y_{T_j} = 1$, and therefore
$L_2 = 1$.  Consequently, on the event $A$, we can see from the definitions
of $(B_i)_{i=0}^{\infty}$ and $(D_i)_{i=0}^{\infty}$ that
$B_M = D_M$ and thus, using part 2 of Lemma \ref{endlem},
$B_K \geq B_M = D_M \geq N$.

Since $(Z_k)_{k=0}^{\infty}$ is a nonnegative supermartingale, it
follows from the Martingale Convergence Theorem (see Corollary 2.11
in chapter 4 of \cite{durr96}) that there exists a random variable $Z$
such that $Z_k \rightarrow Z$ a.s. and $E[Z] \leq E[Z_0]$.
On the event $A$, if $j > K$ then $T_j = T_K$ and thus $Z_j = Z_K = Z$.
Hence, using (\ref{defW}),
\begin{align}
Z1_A &= Z_K 1_A = \alpha^{-B_K} N_K (1-w)^{-\log (X_{T_K} - Y_{T_K})} 1_A
\nonumber \\
&\geq \alpha^{-N}(1-\alpha)(1-w)^{-\log(l/2)} 1_A \geq C^{-1} 1_A. \nonumber
\end{align}
It follows that $E[Z] \geq C^{-1} P(A)$.  Thus,
$Q(f) - Q(f - \epsilon) = P(A) \leq CE[Z] \leq CE[Z_0] =
C (1-w)^{-\log \epsilon}$, as claimed.

\subsection{Proof of part 3 of Proposition \ref{mainprop}}

Let $D^1 = S \cap [0,l]$ and $D^2 = S \cap [1- \ell, 1]$.
Define a sequence of stopping times $(\tau_k)_{k=0}^{\infty}$ by
$\tau_0 = 0$ and $\tau_{k+1} = \inf\{n > \tau_k: X_n \in
[0, \ell] \cup [1- \ell, 1]\}$ for all $k \geq 0$.  Then define
$$D_k = \{f: P(X_{\tau_j} = 1-l \mbox{ for some }j \leq k|X_0 = f) > 0\}.$$
Let $D_k^1 = D_k \cap [0, \ell]$ and $D_k^2 = D_k \cap [1- \ell, \ell]$.
Note that $D^1 = \bigcup_{k=0}^{\infty} D_k^1$ and
$D^2 = \bigcup_{k=0}^{\infty} D_k^2$.  We have
$D_0^1 = \emptyset$ and $D_0^2 = \{1- \ell \}$.  For $k \geq 1$,
\begin{align}
D_k^1 &= \{f \in [0, \ell]: P(X_{\tau_1} \in D_{k-1}|X_0 = f) > 0\} \cup D_{k-1}^1
\nonumber \\
D_k^2 &= \{f \in [1- \ell,1]: P(X_{\tau_1} \in D_{k-1}|X_0 = f) > 0\}
\cup D_{k-1}^2. \nonumber
\end{align}

Suppose $X_0 = f$.  If $f \in (\ell, 1- \ell)$, then $s(f) = \ell$ for
all $k < \tau_1$.  Therefore $X_{\tau_1} = f + n \ell$ for some
$n \in \Z$.  If instead $f \in [0, \ell]$, then $s(f) = f$, in which case
either $X_1 = X_{\tau_1} = 0$ or $X_1 = 2f$.
If $X_1 = 2f$, then $X_{\tau_1} = 2f + n \ell$ for some
$n \in \Z$, where $n = 0$ if $2f \in [0,l] \cup [1- \ell, 1]$.  Likewise,
suppose $f \in [1- \ell, \ell]$.  Then $s(f) = 1-f$, so either
$X_1 = X_{\tau_1} = 1$ or $X_1 = 2f-1$.  If $X_1 = 2f-1$, then
$X_{\tau_1}(f) = 2f - 1 + n \ell$ for some $n \in \Z$, where $n = 0$
if $2f-1 \in [0, \ell] \cup [1- \ell, 1]$.

We claim that if $f \in D^1 \cup D^2$, then there exist integers
$a$, $b$, and $c$ such that $f = 2^{-c}(a + b \ell)$.  Furthermore,
we claim that if $f \neq 1- \ell$, then we can choose $a$, $b$, and $c$
such that $c \geq 1$, $a \geq 1$, $a$ or $b$ is odd, and $a \geq 2$
if $f \in D^2$.  We will prove these claims by induction on $k$.
Note that $D_0 = \{1- \ell \}$, so for $f \in D_0$
we can take $a = 1$, $b = -1$, and $c = 0$.
Now, suppose our claims hold when $f \in D_{k-1}$, where
$k \geq 1$.  To show that our claims hold when $f \in D_k$,
we consider two cases.

First, suppose $f \in D_k^1 \setminus D_{k-1}$.  Then
$P(X_{\tau_1} = g|X_0 = f) > 0$ for some $g \in D_{k-1}$.
Since $0 \notin S$, we must have
$2f + n \ell = g$, or equivalently $f = (g - n \ell)/2$, for some
$n \in \Z$ and $g \in D_{k-1}$.  If $g = 1- \ell$, then
$f = (1 - (n+1)\ell)/2$, so $f = 2^{-c}(a + b \ell)$, where $a = 1$,
$b = -(n+1)$, and $c = 1$.  If $g \neq 1- \ell$, then
$g = 2^{-c}(a + b \ell)$, where $c \geq 1$, $a \geq 1$, and
$a$ or $b$ is odd.  Then $f = 2^{-(c+1)}(a + b \ell - 2^c n \ell) =
2^{-(c+1)}(a + (b - 2^c n)\ell)$.  Note that $c+1 \geq 1$, $a \geq 1$,
and $b - 2^c n$ is odd if $b$ is odd, so $a$ or $b - 2^c n$ is odd.

Next, suppose $f \in D_k^2 \setminus D_{k-1}$.  Then
$P(X_{\tau_1} = g|X_0 = f) > 0$ for some
$g \in D_{k-1}$.  Since $1 \notin S$,
we have $2f - 1 + n \ell = g$, or equivalently $f = (1 + g - n \ell)/2$,
for some $n \in \Z$ and $g \in D_{k-1}$.  If $g = 1- \ell$, then
$f = (2 - (n+1)\ell)/2$.  If $n+1$ were even, then $f = 1- m \ell$ for some
positive integer $m$; since $D^2 \subseteq [1- \ell, 1]$, we would have
$f \in \{1- \ell, 1\}$, which is a contradiction because $1 \notin D_k^2$
and $1- \ell \in D_{k-1}$.  Therefore, $n+1$ is odd, so
$f = 2^{-c}(a + b \ell)$, where $c = 1$, $a = 2$,
and $b$ is odd.  If instead $g \neq 1- \ell$, then
$g = 2^{-c} (a + b \ell)$, where $c \geq 1$, $a \geq 1$, and $a$ or $b$
is odd.  Then $f = 2^{-(c+1)}(2^c + a + b \ell - 2^c n \ell) =
2^{-(c+1)}[(2^c + a) + (b - 2^c n) \ell]$.  Note that $c+1 \geq 1$,
$2^c + a \geq 2$, and either $2^c + a$ or $b - 2^c n$ is odd because
$2^c$ and $2^c n$ are even and either $a$ or $b$ is odd.  It now
follows by induction that our claims hold for all $f \in D^1 \cup D^2$.

Since $\ell < 1/2$, we can choose a positive integer $m$ such that
$1 - m \ell \in (\ell, 2\ell]$.  We can then choose positive integers $d$ and $n$
such that $2^{-d}(1 - m \ell) < 1 - 2 \ell$ and $2^{-d}(1 -  \ell) + n \ell \in
(1 -2 \ell, 1- \ell)$.  Let $f = 2^{-d}(1 - m \ell) + n \ell$.  Note that $f \in (\ell, 1- \ell)$.
Also, $f - \ell = 2^{-d}(1 - m \ell) + (n-1)\ell$.
Suppose a gambler who starts with a fortune of $f - \ell$ loses the
first $n - 1$ bets, then wins the next $d + m - 1$.
After the $n-1$ losses,
the gambler's fortune will be $2^{-d}(1 - m \ell)$.  Then after $d$ wins,
the fortune will be $1-m \ell$.  After $m-1$ additional wins,
the gambler's fortune will be $1- \ell$.  Consequently,
$$P(X_{n+m+d-2} = 1- \ell|X_0 = f - \ell) > 0,$$
which means $f - \ell \in S$.  We now show by contradiction that
$f + \ell \notin S$, which will complete the proof.
Suppose $f + \ell \in S$.  Since $f + \ell > 1- \ell$, there exist
integers $a$, $b$, and $c$ such that $a \geq 2$, $c \geq 1$,
$a$ or $b$ is odd, and $f + \ell = 2^{-c}(a + b \ell)$.  We also have
$f + \ell = 2^{-d}(1 + (2^d(n+1) - m)\ell)$.  Therefore, $2^d(a + b \ell) =
2^c[1 + (2^d(n+1) - m)\ell]$, and so $2^d a - 2^c =
[2^c(2^d(n+1) - m) - 2^d b]\ell$.  Since $\ell$ is irrational, we must have
$2^d a - 2^c = 2^c (2^d (n+1) - m) - 2^d b = 0$.  Thus, $2^d a = 2^c$,
and since $a \geq 2$, it follows that $a$ is even and $c > d$.
Therefore, $b$ is odd and
$b = 2^{c-d}(2^d(n+1) - m)$, which is a contradiction.

\subsection{Obtaining Theorem \ref{mainth} from Proposition \ref{mainprop}}

Suppose $\ell$ is irrational.  By part 3 of Proposition \ref{mainprop}, there
exists $f_0 \in (\ell, 1-\ell)$ such that $f_0 - \ell \in S$ and
$f_0 + \ell \notin S$.  Let $f = f_0 - \epsilon$, where $0 < \epsilon < \ell$ 
and $\epsilon$ is small enough that $f \in (\ell, 1- \ell)$.  We will show that
for sufficiently small $\epsilon$, we have
\begin{equation}
wQ(f + \ell - \epsilon) + (1-w)Q(f - \ell + \epsilon) > Q(f),
\label{fin1}
\end{equation}
which implies Theorem \ref{mainth} because $s(f) = \ell$.  Note that
\begin{equation}
wQ(f + \ell - \epsilon) + (1-w)Q(f - \ell + \epsilon) =
wQ(f_0 + \ell - 2 \epsilon) + (1-w) Q(f_0 - \ell).
\label{fin2}
\end{equation}
Since $f \mapsto Q(f)$ is nondecreasing, we have
\begin{equation}
Q(f) = w Q(f + \ell) + (1-w) Q(f- \ell) \leq w Q(f_0 + \ell) + (1-w)
Q(f_0 - \ell - \epsilon). 
\label{fin3}
\end{equation}
Since $f_0 - \ell \in S$, it follows from part 1 of Proposition \ref{mainprop}
that there exists a constant $C > 0$ such that
\begin{equation}
Q(f_0 - \ell) - Q(f_0 - \ell - \epsilon) \geq C (1-w)^{-\log \epsilon}.
\label{fin4}
\end{equation}
Let $C_0 = C (1-w)$.  Since $f_0 + \ell \notin S$, part 2 of Proposition
\ref{mainprop} implies that for sufficiently small $\epsilon$, we have
\begin{equation}
Q(f_0 + \ell) - Q(f_0 + \ell - 2\epsilon) \leq C_0 (1-w)^{-\log 2\epsilon}
= C(1-w)^{-\log \epsilon}.
\label{fin5}
\end{equation}
Let $B = C(1-w)^{-\log \epsilon}$.  Equations (\ref{fin2})-(\ref{fin5})
imply $$wQ(f + \ell - \epsilon) + (1-w)Q(f - \ell + \epsilon) - Q(f)
\geq -wB + (1-w)B = (1-2w)B > 0$$
for sufficiently small $\epsilon$, which gives (\ref{fin1}).

\bigskip
\begin{center}
{\bf {\Large Acknowledgments}}
\end{center}

The author thanks Lester Dubins for introducing him to this area and for
helpful discussions regarding this work.  He also thanks David Gilat
for bringing to his attention the reference \cite{hps72} and two
anonymous referees for their useful comments.

\end{document}